\documentclass[11pt]{article}

\usepackage{latexsym,amsmath,amsfonts,amssymb,amstext,theorem,euscript,array,enumerate,mathrsfs,color}

\newtheorem{Theorem}{Theorem}[section]

\numberwithin{equation}{section}

\def\Dim{\noindent\hbox{{\bf Proof.}$\;\; $}}
\def\finedim{{\hfill\hbox{\enspace${ \square}$}} \smallskip}

\usepackage[margin=0.7in]{geometry}
\def \E {\mathbb{E}}
\def \P {\mathcal{P}}
\def \e {\varepsilon}

\def\F {\mathcal{F}}
\def \R {\mathbb{R}}
\newtheorem{theorem}{Theorem}[section]
\newtheorem{lemma}[theorem]{Lemma}

\newtheorem{hypothesis}[theorem]{Hypothesis}

\newtheorem{remark}[theorem]{Remark}
\numberwithin{equation}{section}
\def\Dim{\noindent\hbox{{\bf Proof.}$\;\; $}}          
\def\finedim{{\hfill\hbox{\enspace${ \square}$}} \smallskip}    
\def\sqr#1#2{{\vcenter{\vbox{\hrule height .#2pt
     \hbox{\vrule width .#2pt height#1pt \kern#1pt \vrule
     width .#2pt} \hrule height .#2pt}}}}
\def\square{\mathchoice\sqr54\sqr54\sqr{4.1}3\sqr{3.5}3}

\begin{document}
\thispagestyle{empty}
\parindent=0pt

\title{Singular Limit of Two Scale Stochastic Optimal Control Problems in Infinite Dimensions by Vanishing Noise Regularization}
\author{Giuseppina Guatteri, \\
Dipartimento di Matematica, \\
Politecnico di Milano, \\
Piazza Leonardo da Vinci 32, \\
20133 Milano,
Italia. \\
 {\tt e-mail: giuseppina.guatteri@polimi.it} \\ \\
Gianmario Tessitore,\\
Dipartimento di Matematica e Applicazioni \\
Universit\`{a} Milano-Bicocca, \\
 via R. Cozzi 53 - Edificio U5,\\ 20125 Milano, Italia.\\
{\tt e-mail: gianmario.tessitore@unimib.it
 }} 

\date{}
\maketitle
\begin{abstract}
In this paper we  study  the limit of the value function for a two-scale,  infinite-dimensional, stochastic controlled system with cylindrical noise and possibly degenerate diffusion. 
The limit  is represented as the value function of a new \textit{reduced} control problem (on a reduced state space).
The presence of a cylindrical noise prevents representation of the limit by viscosity solutions of HJB equations as in \cite{Swiech2021} while degeneracy of diffusion coefficients  prevents representation as a classical BSDE as in \cite{GuaTess2019AMO}. We use a  ˝vanishing noise" regularization technique.
\end{abstract}

\section{Introduction}

This  paper studies the limit of the value functions for a sequence of optimal control problems with state equation represented by  the following system of stochastic differential equations
\begin{equation}
\label{lentovelocecontrollata-intro}
\begin{cases}
dX_t= AX_t \, dt  +  b(X_t,Q_t,u_t)dt \,+ R (X_t) d{W}^1_t, & X_0=x_0, \\ 
 \e dQ_t= (BQ_t+ F(X_t,Q_t)  +  G\rho(u_t))\;dt+ \e^{1/2} G \, d {W}^2_t, & Q_0=q_0.
\end{cases}
\end{equation}
and cost represented by the following functional:
$$J^{\epsilon}(x_0,q_0,u)=\mathbb{E}\left[\int_0^1 l(X_t,Q_t,u_t)dt+h(X_1)\right].$$
We notice that the small constant $\e$ in the second equation modelizes the fact that $Q$ evolves  quicker than $X$, with a ratio $\frac{1}{\e}$ between the two velocities. Our goal is to represent the limit of the value functions of these problems as the speed factor diverges.

In \eqref{lentovelocecontrollata-intro} both $X$ and $Q$  take values in an Hilbert space. Moreover $A$ and $B$ are unbounded linear operators, $u$ represents the control, $  ({W}^1_t)_{t\geq 0}$,  $  ({W}^2_t)_{t\geq 0}$ are infinite dimensional independent cylindrical Wiener processes, $b$, $F$, $\rho$ are $R$ are functions  and $G$ is a  bounded linear operator satisfying suitable assumptions including dissipativity of $B+F(x,\cdot)$.
The main feature of this paper is that we allow both $W^1 $ and $W^2$ to be cylindrical while we do not require any regularizing or nondegeneracy assumptions to hold on $G$ or $R(\cdot)$. We mention here that throughout the paper the control problems are formulated in a weak form particularly suitable   to be studied by Backward Stochastic Differential Equations (BSDEs for short).

\medskip 

Several papers have been devoted to the characterization of  limits  of singular stochastic control problems in finite dimensional spaces. Beside the pioneering results based on direct computations in specific situations (see, for instance, \cite{KabPer} and  \cite{KabPer1}) the general approach in finite dimensional cases (see   \cite{AB3}, \cite{AB2}, \cite{AB1}, \cite{AC1}, \cite{AC2})  relies on the representation of the value function  as  \textit{viscosity} solution of a suitable HJB equation and  on a  convergence result, as $\epsilon \rightarrow 0$, for  the solution of such HJB equations  to a  \textit{reduced} nonlinear parabolic equation. The value of the (viscosity) solution of the limit PDE is then the desired limit. The well known technical difficulties  nested in the proof of comparison principle for solution of infinite dimensional HJB equations prevents a direct extension of the previous results to the case of Hilbert-valued,  two scale, controlled stochastic systems. 

At our best knowledge the first paper to address the  problem in an infinite dimensional framework is \cite{GuaTess2019AMO} where the value function is represented through a BSDE and the result is obtained by convergence of a class of singularly perturbed BSDEs. The main limitation of the results in \cite{GuaTess2019AMO}  is that we need to assume non degeneracy of the noise in the slowly evolving equation (in \cite{GuaTess2019AMO} $R$ is indeed independent of $x$ and, more important, is invertible). Then in \cite{Swiech2021} the viscosity solution approach is adapted to the Hilbert space case by a deep analysis of the  necessary  technical assumptions. A rather general class of two scales systems can be considered in this last paper with the only remaining obstruction on the covariance of the noises that must be of finite trace. 

As we have already mentioned here both $(W^1)$ and $(W^2)$ are cilindrical and we do not assume non degeneracy nor on $R$ nor on $G$ (nevertheless in the equation for $Q$ the \textit{structure condition}, allowing to apply Girsanov transform, has to hold). Since the two previously mentioned possible representations of our singular limit (the one through a BSDE and the one through a viscosity solution of a HJB equation) seem not to be available in the present case we try to represent it as the value function of a \textit{reduced} control problem. By the way we notice that such a representation somehow lies underneath both the above mentioned ones. It is also  worth mentioning that any notion of solution of HJB equation stronger than viscosity, and consequently any  direct representation of the limit value function by a standard BSDE, appears here to be excluded by the lack of regularity that one has to expect for the value functions of degenerate stochastic control problem.

To compare the class of state equations that fall into the framework of this paper with the ones treatable in \cite{GuaTess2019AMO} and \cite{Swiech2021} let us consider the following two scale system of controlled reaction diffusion  SPDEs in one space dimension driven by space time noises. We refer, for instance to \cite{DpZ1} Section 11.2 for the abstract formulation and precise assumptions on the coefficients. We just have to mention that  $m$ is a positive constant and the Lipschitz constant of 
$f$  with respect to $\mathcal{Q}$ is smaller then $m$. 

\begin{equation}\label{sistemaesempio}\left\{\begin{array}{l}
\displaystyle\frac{\partial}{\partial t} \mathcal{X}^\e(t,x) = \frac{\partial^2}{\partial x^2} \mathcal{X}^\e(t,x) +
 b(\mathcal{X}^\e(t,x),\mathcal{Q}^\e(t,x), u(t,x)) +
 \sigma(x,\mathcal{X}^\e(t,x) )  \frac{\partial}{\partial t} {\mathcal{W}}^1(t,x),   \\ \\\displaystyle
\e\frac{\partial}{\partial t} \mathcal{Q}^\e(t,x)\! =\!(\frac{\partial^2}{\partial x^2}- m) \mathcal{Q}^\e(t,x)\! + \!
f(\mathcal{X}^\e(t,x),\mathcal{Q}^\e(t,x))+  \rho(x) r(u(t,x)) 
+ \e ^{1/2}\rho(x) \frac{\partial}{\partial t}
{\mathcal{W}}^2(t,x),   \\ \\
\mathcal{X}^\e(t,0)=\mathcal{X}^\e(t,1)=\mathcal{Q}^\e(t,0)=\mathcal{Q}^\e(t,1)=0, \,\\ \\\mathcal{X}^\e(0,x)= \mathcal{X}^0(x), \ \mathcal{Q}^\e(0,x)= \mathcal{Q}^0(x),\qquad\qquad\qquad\qquad  t \in [0,1], \ x \in [0,1],
\end{array}\right.
\end{equation}
To fit the assumptions in \cite{GuaTess2019AMO} one should assume $\sigma $ to be independent of $\mathcal{X}$ and bounded away from 0 while to fit the assumptions in \cite{Swiech2021} one should assume $({\mathcal{W}}^i)$, $i=1,2$ to be colored in space. Here we can take a general $\sigma$ (possibly vanishing ) and 
space-time white noises $({\mathcal{W}}^1)$, $i=1,2$. 

The approach we present here is to regularize equation \eqref{lentovelocecontrollata-intro} adding an extra noise with small parameter in the equation for $X$. We obtain a non degenerate  singular  control problem (see equation \eqref{lentovelocecontrollataviscosa})  that can be treated using the results in \cite{GuaTess2019AMO}.
The point is that, in this way, we have two handle  a system depending on two parameters (the original speed ratio $\e$ and the new \textit{small noise} parameter). We show that we can interchange the limits and let first $\e \rightarrow 0$ for a fixed small noise parameter. Proceeding like that (adapting the arguments in \cite{GuaTess2019AMO}) we end up with a forward backward system depending on the small parameter $\eta$.
\begin{equation}\label{forbacetaintro}
\begin{cases}
d X_t =   AX_t \, dt + R(X_t)\, d {W}^1 _t + \eta \, d {B}_t, \qquad t \in [s,1], \\
-d Y _t  = \lambda({X}_t,\eta^{-1}{Z}_t)\,dt - \, {Z}_t\, d {W}^1_t + \, {Z}_t\, d {B}^1_t,   \\
X_s= x, \quad {Y} _1=h(X_1).
\end{cases} 
\end{equation}
We notice that in the above equation the nonlinearity $\lambda$ is itself the value function of a suitable ergodic optimal control problem, roughly speaking, for the second equation in \eqref{lentovelocecontrollata-intro}  with frozen $X$.

The idea is then to see, by Fenchel duality, the $Y$ in equation  \eqref{forbacetaintro} as the value function of an auxiliary control problem. The key issue, at this level, is to construct this new control problem in such a way that its running cost has enough regularity to allow the final passage to the limit as the small noise regularization vanishes and, eventually, to get the main result of this paper (see Theorem \ref{main-caratt}). In addition the value function of our reduced control problem can be  shown to coincide with the minimal solution of a Backward Stochastic Differential Equation with constraints on the martingale term (see Remark \ref{constrainedBSDEs} here and \cite{KharroubiPham}, \cite{CossoGuatteriTessitore}).

The paper is organized as follows: in Section 2 we introduce general notation, in Section 3 we formulate the control problems introducing the weak formulation that will be used throughout the paper, in Section 4 we introduce the small noise  regularization of the system, in Section 5 we prove that we can change order between the limit with respect to the speed ratio parameter $\e$ and the small noise parameter $\eta$, finally in Section 6 we prove our mail result.

\section{Notation }
Given a Banach space $E$,
the norm of its elements $x$ will be denoted
by $|x |_E$, or even by $|x|$ when no confusion is possible. If $F$ is
another Banach space, $L(E,F)$ denotes the space of bounded linear
operators from $E$ to $F$, endowed with the usual operator norm. When $F=\mathbb{R}$ the dual space $L(E,\mathbb{R})$ will be denoted by $E^*$.  
The letters $\Xi$, $H$ and $K$ will always be used to denote Hilbert spaces.
The scalar product is denoted $\langle \cdot, \cdot \rangle$, equipped with a
subscript to specify the space, if necessary. All  Hilbert
spaces are assumed to be real and separable and  the dual of a Hilbert space will never be identified with the space itself. By $L_2(\Xi,H)$ and
$L_2(\Xi,K)$
 we denote the spaces of Hilbert-Schmidt operators from $\Xi $ to $H$
 and to $K$, respectively. Finally  $\mathcal{G}(K,H)$ is the space of all  Gateaux differentiable mappings $\phi$
 from $K$ to $H$ such that the map $(k,v)\rightarrow \nabla \phi(k)v$ is continuous from $K\times K$ to $H$; see \cite{FuTessAOP2002} for details.
 \vspace{5pt}



Next we define  the following classes of stochastic processes with values in
a Hilbert space $V$.  Given an arbitrary time horizon  $T$, constant $p\geq 1$ and a generic filtered space $(\Omega^0,\mathcal{E}^0,(\mathcal{F}^0_t)_{t\in [0,T]},\mathbb{P}^0)$:
\begin{itemize}
\item $L^p_{\mathcal{\F}^0} (\Omega^0\times [0,T];V)$ denotes the space of
equivalence classes of processes $Y \in L^p (\Omega\times
[0,T];V)$ admitting a predictable version. It is endowed with the norm
\[ |Y|_p= \Big(\E^0 \int_0^T |Y_s|^p \, ds\Big)^{1/p}. \]
\item $L^p_{\mathcal{\F}^0}(\Omega^0;C([0,T];V))$
     denotes the space of
    adapted processes $Y$ with continuous paths in $V$, such
    that the norm
    \[  \|Y\|_p  = (\E^0 \sup _{s \in [0,T]} |Y_s|^p)^{1/p}\]
    is finite. The elements of $L^p_{\mathcal{\F}^0}(\Omega^0;C([0,T];V))$
    are identified up to indistinguishability.
\end{itemize}


\section{Setting of the problem and statement of the main result}
\label{setting}

Let $H$, $K$ and $\Xi$ separable Hilbert spaces and $U$ a separable
metric space. 
We denote  by $\mathbb{S} ^{1,2}$ ($\mathbb{S}$ stands for \textit{Setting}) the class of all {$6$-uples $\mathbb{U}=(\Omega, (\F_t), \mathbb{P}, (W^1_t), (W^2_t), (u_t))$},  where $(\Omega, \F, (\F_t))$ is a filtered complete probability space, $(W^1_t), (W^2_t)$ are two independent, $\Xi$-valued,  $(\F_t)$-Wiener processes and $u$ is an $(\F_t)$ predictable process taking values in $U$.  When needed, we will  add the mark  $\mathbb{U}$ to each term to avoid confusion.

Given $x_0 \in H, q_0 \in K, \e >0$, and $\mathbb{U} \in \mathbb {S} ^{1,2}$,  we consider the following two scale state equation in $H \times K$:
\begin{equation}
\label{lentovelocecontrollata}
\begin{cases}
dX_t= AX_t \, dt  +  b(X_t,Q_t,u_t)dt \,+ R (X_t) d{W}^1_t, & X_0=x_0, \\ \\
 \e dQ_t= (BQ_t+ F(X_t,Q_t)  +  G\rho(u_t))\;dt+ \e^{1/2} G \, d {W}^2_t, & Q_0=q_0.
\end{cases}
\end{equation}
that has, under Hypothesis 3.1-3.6 listed below,  a unique mild solution belonging to $L^p_{\mathcal{\F}^\mathbb{U}}(\Omega^{\mathbb{U}};C([0,T];H))$, $p\geq 1$ that we denote by $ X^{\e, \mathbb{U}}$, see \cite[Lemma 3.9 and Lemma 3.10]{GuaTess2019AMO}.  We omit reference to initial state $(x_0,q_0)$ trying to ease the notation (when we will need to show such dependence we will explicitly mention it).
We introduce the following cost functional to minimize
\begin{equation} \label{cost_epsilon}
J^{\e}(x_0,q_0,\mathbb{U})= \E^\mathbb{U}\left[\int_0^1 l(X^{\e,\mathbb{U}}_t, Q^{\e,\mathbb{U}}_t, u^\mathbb{U}_t)dt + h (X^{\e,\mathbb{U}}_1)\right],
\end{equation}
where $\E ^{\mathbb{U}}$ denotes the expectation with respect to the probability $\mathbb{P}$ in $\mathbb{U}$.
\medskip 

We  make the following general assumptions fixing, in the mean time, three constants $M>0$, $L>0$ and $\gamma\in [0,1/2)$ that will not be changed throughout the paper.
\begin{hypothesis}\label{A.1}
$A: D(A)\subset H \to H$
is a 
linear, unbounded  operator that generates a  $C_0$- semigroup $\{
e^{tA} \}_{t \geq 0}$, such that  $|e^{tA}|_{L(H,H)} \le M _A e ^{ \omega_A t},  t \geq 0$ for some positive constants $M_A$ and $\omega_A$.
$B: D(B)\subset K \to K$ is a 
linear, unbounded  operator that generates a  $C_0$- semigroup $\{ e^{tB} \}_{t \ge 0}$ such that $|e^{tB}|_{L(K,K)} \le   M_B e ^{ \omega_B t}, t \geq 0$
 for some $ M_B, \omega_B>0$. \\
Moreover there exist $C>0$ s.t.:
\begin{align}
&|e^{sA}|_{L_2(H,H)}+ |e^{sB}|_{L_2(K,K)}\nonumber \leq C
s^{-\gamma},\quad \forall s\in [0,1].\nonumber 
\end{align}
\end{hypothesis}
\begin{hypothesis}\label{A.2}
 The functions $b: H
\times K \times U  \to H $ and $F:  H
\times K \to K$ are measurable and:
\begin{equation*}
   |b(x,q,u)|\leq M,\quad   |b(x,q,u)- b(x',q',u)| \leq L (|x-x'|+|q-q'|),\qquad \qquad \forall \, q,q' \in K,  x, x' \in H, \, u \in U,
\end{equation*}
\begin{equation*}
| F(x,q)-
F(x',q')| _K \leq L( |x-x'|_H + |q-q'|_K)\qquad \qquad \forall \, q,q' \in K,  x, x' \in H.
\end{equation*}
Moreover we assume that,  $F (x,\cdot)$ is Gateaux differentiable, more precisely, $F (x,\cdot) \in \mathcal{G}^1(K,K)$, $\forall x\in H$.
\end{hypothesis}
\begin{hypothesis}\label{A.3} $B +F$ is  dissipative i.e.
there exists some $\mu >0$ such that:
\[ \langle B q+F (x,q) - (Bq' +F(x,q')), q-q' \rangle
 \leq -\mu |q-q'|^2, \]
for all $
x\in H,
 q,q' \in D(B)$.
\end{hypothesis}
\begin{hypothesis} \label{A.5}
$R\colon H \rightarrow  L(\Xi,H)$  is a bounded Lipschitz map. Moreover for all  $x,x'\in H$:
$$
\big|e^{sA}R (x)-e^{sA}R (x')\big|_{L_2(\Xi,H)} \ \leq \ \frac{L }{ s^\gamma}  |x-x'|_{H},\qquad 
\hbox{
for all $s\in (0,1)$.}$$
\end{hypothesis}

\begin{hypothesis}\label{A.6}
$G \in {L}(\Xi;K)$.
\end{hypothesis}

\begin{hypothesis}\label{A.7}
 The functions $l: H\times K\times U\rightarrow \mathbb{R}$ and $h: H \rightarrow \mathbb{R}$ are measurable and  satisfy the assumptions below, moreover
$$ |l(x,q,u)- l(x',q',u)| \leq L (|x-x'|+|q-q'|),\qquad \qquad \forall \, q,q' \in K,  x, x' \in H, \, u \in U,$$
$$ |h(x)- h(x')| \leq L |x-x'|,\qquad \qquad \forall x, x' \in H, $$
$$ |l(x,q,u)|, |\rho(u)|, |h(x)| \leq M,\qquad \qquad  \forall q \in K, x \in H, \, u \in U. $$
\end{hypothesis}

We are interested in studying the limit of the value function $ V^\e (x_0,q_0) $, 
\begin{equation}\label{valuefunction_orig}
{V}^\e (x_0,q_0)=:\inf_{\mathbb{U} \in \mathbb{S}^{1,2}} J^\e (x_0,q_0,\mathbb{U})
\end{equation} 
as the ratio $\e $ between the speed of the slow component and the speed of the fast one tends to $0$. Namely  we shall provide a representation of this limit by a \textit{ reduced } stochastic control problem.


\section{Small noise \label{sec-statement} approximations of the two scale problem}

In order to regularize our initial problem we introduce a {\em vanishing}  noisy term  in \eqref{lentovelocecontrollata}.  To do that we have to modify our class of settings.

Namely  we denote  by $\mathbb{S} ^{1,2,B}$ the class of $7$-uples  $\mathbb{U}^B=(\Omega, (\F_t), \mathbb{P}, (W^1_t), (W^2_t), (B_t) ,(u_t))$, 
where, beside the forementioned elements, there is a third  $(\F_t)$-Wiener process  $(B_t)$, independent of $ (W^1_t,W^2_t)$.

Then given $x_0,q_0$ and a setting $\mathbb{U}^{B} \in \mathbb{S}^{1,2,B}$, for every $\eta \geq 0$, let us consider the following regularized two scale state equation:
\begin{equation}
\label{lentovelocecontrollataviscosa}
\begin{cases}
dX_t= AX_t \, dt  +  b(X_t,Q_t,u_t)dt \,+ R (X_t) dW^1_t + \eta\,  d B_t, & X_0=x_0, \\ \\
 \e dQ_t= (BQ_t+ F(X_t,Q_t)  +  G\rho(u_t))\; dt+ \e^{1/2} G \, d W^2_t, & Q_0=q_0.
\end{cases}
\end{equation}

Such system has a unique mild solution, indeed following \cite[Lemma 3.9 and Lemma 3.10]{GuaTess2019AMO} we have that  for every $\e >0, \eta >0$  there exists a unique  couple of processes $(X ^{\e,\eta,\mathbb{U}^B}, Q^{\e, \eta, \mathbb{U}^B})$,  with $ X  ^{\e,\eta,\mathbb{U}^B}$  belonging to $ L^p_{\mathcal{F}^{\mathbb{U}^B}}(\Omega^{\mathbb{U}^B};C([0,1];H))$ and $Q^{\e, \eta, \mathbb{U}^B}$ in $L^p_{\mathcal{F}^{\mathbb{U}^B}}(\Omega^{\mathbb{U}^B};C([0,1];K))$ .

Moreover, see  again  \cite[Lemma 3.9 and Lemma 3.10]{GuaTess2019AMO}, also see the proof of Theorem \ref{limite epsilon fisso}  below, the following estimates hold:
\begin{equation}\label{stimaslow}
 \E^{\mathbb{U}^B}(\sup_{t\in [0,1]}|X^{\e,\eta,\mathbb{U}^B}_t|^p)\leq c_p(1+|x_0|^p), \qquad x_0 \in H, \qquad \forall p \geq 1
 \end{equation}
\begin{equation}
\sup_{t\in [0,1]} \E^{\mathbb{U}^B} |{Q}^{\e,\eta,\mathbb{U}^B}_s|^p \leq k_p(1+ |q_0|^p), \qquad q_0 \in K, \qquad \forall p \geq 1.
\end{equation}
with constant $ c_p$ and $ k_p$ independent from $\e$ and $\eta$.
\medskip

We also consider the analogue of our control problem in this enriched and regularized situation  and the  corresponding value  function $V^{\e,\eta}$. Namely:  
\begin{equation} \label{cost_epsilon_eta}
\begin{split}
& J^{\e,\eta}(x_0,q_0,\mathbb{U}^B)= \E^{\mathbb{U}^B}\left[\int_0^1 l(X^{\e,\eta,\mathbb{U}^B}_t, Q^{\e,\eta,\mathbb{U}^B}_t, u^{\mathbb{U}^B}_t)dt + h (X^{\e,\eta, \mathbb{U}^B}_1)\right].\\
& V^{\e,\eta}(x_0,q_0)=\inf_{\mathbb{U}^B \in \mathbb{S} ^{1,2,B}} J^{\e,\eta}(x_0,q_0,\mathbb{U}^B)
\end{split}
\end{equation}
It is straightforward to verify that, fixed $(x_0,p_0)$, the set $\{V^{\e,\eta}(x_0,q_0) : \e>0, \eta >0 \}  $ is bounded. 

\begin{remark}\label{CommentoUeUB}

Given a setting $\mathbb{U}^B$ in  $\mathbb{S} ^{1,2,B}$ we define  the setting $P\mathbb{U}^B$  in $\mathbb{S} ^{1,2}$ as the  setting obtained by omitting the process $B$. Namely $$ P (\Omega, (\F_t), \mathbb{P}, (W^1_t), (W^2_t), (B_t),  (u_t))  = (\Omega, (\F_t), \mathbb{P}, (W^1_t), (W^2_t),  (u_t)) $$
In particular the control $u$ is the same in $\mathbb{U}^B$ and $P\mathbb{U}^B$, moreover $
(X^{\e,0,\mathbb{U}^B}, Q^{\e,0,\mathbb{U}^B})\equiv
(X^{\e,P\mathbb{U}^B}, Q^{\e,P\mathbb{U}^B}).
$ On the other hand given $\mathbb{U} \in \mathbb{S} ^{1,2}$ we define by $R \mathbb{U} \subset \mathbb{S} ^{1,2,B}$  the set, always non empty, of all settings $\mathbb{U} ^B$ obtained from $\mathbb{U}$, by choosing  any $(\Omega ^0, \F ^0, \mathbb{P}^0, (\F^0_t), B)$, where  $(\Omega ^0, \F ^0, \mathbb{P}^0, (\F^0_t))$, is a filtered complete probability space and  $(B)$ is an  $(\F^0_t)$-Wiener process, and setting 
$$ U^B =  (\Omega \times \Omega ^0, (\F \otimes \F^0),  (\F_t \otimes \F^0_t), \mathbb{P} \otimes \mathbb{P}^0, (W'^1_t), (W'^2_t),  (B'_t),(u'_t) ).$$
In the above formula ${W}'^i(\omega,\omega^0):= W^i(\omega)$, $i=1,2$; $u'(\omega,\omega^0)= u(\omega)$, $B'(\omega,\omega^0):= B(\omega^0)$, for every $(\omega,\omega^0)\in\Omega \times \Omega^0$. We  notice that if $ \mathbb{U}^B \in R \mathbb{U}$ then the law of $(X^{\e, \mathbb{U}}, Q^{\e, \mathbb{U}},u  ^ \mathbb{U})$ under $ \mathbb{P} ^{\mathbb{U}}$ coincides with the law of  $(X^{\e,0, \mathbb{U}^B}, Q^{\e, 0,\mathbb{U} ^B},u^{\mathbb{U} ^B})$ under $\mathbb{P}^{\mathbb{U}^B} $.

Thus 
$J^{\e, 0} (x_0,q_0,\mathbb{U} ^B)= J^{\e} (x_0,q_0,\mathbb{U}) $, for every $ \mathbb{U} \in \mathbb{S} ^{1,2}$ and   $ \mathbb{U}^B  \in R \mathbb{U} \subset \mathbb{S} ^{1,2,B} $.
\end{remark}
The above remark entitles us to replace our original control problem with the enriched one in the trivial case $\eta=0$. Namely we have the not very  surprising equality: 
\begin{lemma}\label{identita value function}
$$V^{\e,0}(x_0,q_0)={V}^{\e}(x_0,q_0), \qquad \forall x_0\in H,\; q_0\in K.$$
\end{lemma}
\Dim 
Thanks to the previous remark, for every $\mathbb{U}^B \in \mathbb{S} ^{1,2,B}$
we have that $J^{\e}(x_0,q_0,P \mathbb{U}^B)= J^{\e,0}(x_0,q_0,\mathbb{U}^B)$.
Conversely if $ \mathbb{U} \in \mathbb{S}^{1,2}$ and $ \mathbb{U}^B \in R \mathbb{U}$ we again have that $J^{\e}(x_0,q_0,\mathbb{U})=J^{\e}(x_0,q_0,P\mathbb{U}^B)= J^{\e,0}(x_0,q_0,\mathbb{U}^B)$. Thus the sets 
$\{ J^{\e}(x_0,q_0,\mathbb{U}): \mathbb{U} \in   \mathbb{S} ^{1,2} \}$ and  $\{ J^{\e, 0}(x_0,q_0,\mathbb{U}^B): 
\mathbb{U}^B \in   \mathbb{S} ^{1,2,B} \}$ coincide and obviously the same is true for their infimum.
\finedim
\smallskip 

Now we fix a reference setting $ (\bar{\Omega}, \bar{\F}, \bar{\mathbb{P}}, (\bar{\F_t}), \bar{W}^1, \bar{W}^2, \bar{B})$, we denote with $ \bar{\F} ^{1,2,B}_t,\bar{ \F} ^{1,B}_t, \bar{\F }^2_t$ the natural filtration generated, respectively, by the processes  $ (\bar{W}^1, \bar{W}^2, \bar{B}), (\bar{W}^1, \bar{B})$,  and $\bar{W}^2$ (always completed with the $ \bar{\mathbb{P}}$ -negligible sets in $\bar{\mathcal{F}}$), eventually by $\bar {\E}$ we will denote the expectation with respect to $\bar{\mathbb{P}}$.

 
We then consider the following system:
\begin{equation}\label{lentoveloceottimo}
\begin{cases}
dX_t= AX_t\,dt  +  R(X_t)\,  d\bar{W}^1_t + \eta\,  d\bar{B}_t, & X_0=x_0, \\ \\
 \e dQ_t= (BQ_t+ F(X_t,Q_t) dt +   \e^{1/2} G \, d\bar{W}^2_t, & Q^\e_0=q_0.
\end{cases}
\end{equation}
That has a unique  mild solution $(X^{\eta},Q^{\e,\eta})$ with $ X^\eta\!\in L^p_{\bar {\mathcal{F}}^{1,B}}(\Omega;C([0,1];H))$, $Q^{\e,\eta} \!\in L^p_{\bar{\mathcal{F}}^{1,2,B}}(\Omega;C([0,1];K))$, $p\geq 1$.

We introduce here the BSDE:
\begin{equation}\label{BSDE-prima}
\begin{cases}
 -d Y_t= \psi^{\epsilon,\eta}(X^\eta_t, Q^{\e,\eta}_t,   {Z}^1_t, {Z}^2_t, \Xi_t ) dt-  {Z}^1_t   d\bar{W}^1_t  -  {Z}^2_t   d\bar{B}_t-   {\Xi}_t   d \bar{W}^2_t, \\ \quad Y _1 \  = h(X_1).
\end{cases}
\end{equation}
with\begin{align*}
    & \psi^{\e,\eta}(x, q,   z_2,v) :  =\inf_{\{ u \in{U}\}} \{ l(x,q,u) + \frac{z_2 b(x,q,u)}{\eta} + \frac{ v}{\sqrt{\varepsilon}} \rho(u)\} =   \psi(x, q, \frac{z_2}{\eta}, \frac{ v}{\sqrt{\varepsilon}}), 
\end{align*}
where $\displaystyle \psi(x,q,  z_2,v)=\inf_{ u \in{U}} \{ l(x,q,u) + z_2 b(x,q,u) +  v \rho(u)\}$ for every  $ x, z_2 \in H, q\in K, v \in \Xi$.

We notice that $\psi$ is Lipschitz in $z_2$ and $v$ uniformly with respect to $x$ and $q$, moreover it is Lipschitz with respect to $x$ and $q$ with Lipschitz constant linearly growing in $|z_2|$. 

By standard BSDE theory (see for instance \cite{FuTessAOP2002}) equation \eqref{BSDE-prima} has a unique solution
 $(Y^{\e, \eta}, Z^{1, \e, \eta}, Z^{2, \e, \eta}, \Xi^{\e, \eta})$ with   $Y^{\e,\eta}\! \in L^2_{\bar{\mathcal{F}}^{1,2,B}}(\Omega;C([0,1];\R))$, $Z^{1,\e,\eta} \!\in L^2_{\bar{\mathcal{F}}^{1,2,B}} (\Omega\times [0,1];\Xi^*)$, $Z^{2,\e,\eta} \!\in L^2_{\bar{\mathcal{F}}^{1,2,B}} (\Omega\times [0,1];H^*)$ and  $\Xi^{\e,\eta} \!\in L^2_{\bar{\mathcal{F}}^{1,2,B}}(\Omega\times [0,1];\Xi^*)$.
 
\medskip

We introduce the space of $U$-valued processes $\mathcal{U}^{1,2,B}$ being progressively measurable w.r.t. $\bar{\F}_{t}^{1,2,B}$. 

The following identification  is a standard result in BSDE theory we report the proof in order to take into account the two different (weak) formulations of the control problem that will be needed below.
\begin{lemma}\label{lemma 4.3}
It holds that:
\begin{equation}
    V^{\e, \eta} (x_0,q_0)= \inf_{u \in \mathcal{U}^{1,2,B}} \E ^{u} \Big( \int_0^1 l (X^\eta_t, Q^{\e,\eta}_t))\, dt + h(X^\eta_1)\Big)= Y^{\e,\eta}_0.
\end{equation}
where $\E ^u$ denotes expectation with respect to the probability $\mathbb{P}^u$ on  $ (\bar{\Omega}, \bar{\F})$ under which
$$ (\bar{W}^1_t,- \int_0^1 \e ^{-1/2}\rho(X^\eta_t, Q^{\e,\eta}_t)\, dt  + \bar{W}^2_t, - \int_0^1  \eta ^{-1 }b(X^\eta_t, Q^{\e,\eta}_t, u_t)\, dt + \bar{B}_t)$$
is a Wiener process on $\Xi \times \Xi  \times H$.
\end{lemma}
\Dim
We begin noticing that, given $u \in \mathcal{U}^{1,2,B}$, we can, starting from the reference setting, build the following setting $\bar{\mathbb{U}}^B\in  \mathbb{S}^{1,2,B}$ by: $$\bar{\mathbb{U}}^B:= (\bar{\Omega}, \bar{\F}, {\mathbb{P}^u}, ({\bar{\F}^{1,2,B}_t}),\bar{W}^1_t,- \int_0^1 \e^{-1/2} \rho(u_t)\, dt  + \bar{W}^2_t, - \int_0^1 \eta ^{-1} b(X^\eta_t, Q^{\e,\eta}_t, u_t)\, dt + \bar{B}_t, u_t).$$
$\bar{\mathbb{U}}^B$ has been constucted in such a way that the corresponding solution to
\eqref{lentovelocecontrollataviscosa} $ X^{\e, \eta,{\mathbb{U}}^B}$ coincides with the solution $X^{ \eta}$  of \eqref{lentoveloceottimo}.
We then have that $$J^{\e,\eta}(x_0, q_0, \bar{\mathbb{U}}^B) = \E^{u}  \Big(\int_0^1 l (X^{\eta}_s, Q^{\e,\eta}_s)\, ds + h(X^{\eta}_1)\Big).$$
So, by definition, we have that:
\begin{equation} \label{versofacile}
    V^{\e,\eta}(x_0,q_0)\leq \inf_{u \in \mathcal{U}^{1,2,B}} \E ^{u} \Big( \int_0^1 l (X^\eta_t, Q^{\e,\eta}_t)\, dt + h(X^\eta_1)\Big).
\end{equation}
To prove the converse we  add and subtract the terms $\displaystyle \int_0^1 l(X^\eta_t, Q^{\e,\eta}_t, u_t)\, dt  $ to system \eqref{lentoveloceottimo}  and compute the mean value with respect to $\mathbb{P} ^u$ to get the usual fundamental relation:
\begin{align*}
 Y^{\e, \eta}_0 &= \E ^u \Big[ 
 \int_0^1 (\psi ^{\e,\eta}( X^\eta_t, Q^{\e,\eta}_s, Z^{1,\e, \eta}_t, Z^{2,\e, \eta}_t, \Xi^{\e, \eta}_t) - l(X^\eta_t, Q^{\e,\eta}_t, u_t)   -\frac{1}{\eta} Z^{2,\e, \eta}_t b(X^{\eta}_t, Q^{\e,\eta}_t,u_t) - \frac{1}{\sqrt{\e}}\Xi ^{\e, \eta}_t \rho(u_t))\, dt   \Big ] \\ &
 + \E ^u \Big[\int_0^1 l (X^{\eta}_t, Q^{\e,\eta}_t,u_t)\, dt + h(X^{\eta}_1)\Big]
\end{align*}
Choosing, in a measurable way, a minimizing  sequence $u^n$  such that 
$$|\psi ^{\e,\eta}( X^\eta_t, Q^{\e,\eta}_t, Z^{1,\e, \eta}_t, Z^{2,\e, \eta}_t, \Xi^{\e, \eta}_t) - l(X^\eta_t, Q^{\e,\eta}_t, u^n_t)   - \frac{1}{\eta}Z^{2,\e, \eta}_t b(X^{\eta}_t, Q^{\e,\eta}_t,u^n_t) - \frac{1}{\sqrt{\e}} \Xi ^{\e, \eta}_t \rho(u^n_t))|\displaystyle \leq \frac{1}{n}$$ (notice that $ (X^\eta , Q^{\e,\eta})$ do not depend on $u$)  we end up, see  for instance \cite{FuTessAOP2002}, with:
\begin{equation}
   Y^{\e, \eta}_0 = \inf_{u \in \mathcal{U}^{1,2,B}} \E ^{u} \Big( \int_0^1 l (X^\eta_t, Q^{\e,\eta}_t))\, dt + h(X^\eta_1)\Big)
\end{equation}
To complete the proof it is enough to show that:
\begin{equation*}
    J^{\e, \eta}(x_0,q_0, \mathbb{U}^B) \geq Y^{\e, \eta }_0, \quad \text{for all}
    \   \mathbb{U}^B \in \mathbb{S} ^{1,2,B}.
\end{equation*}
The key point is the well known observation that the law of the solution to the forward backward system \eqref{lentoveloceottimo} does not depend on the specific setting (the solution is obtained by a Picard iteration argument that conserves the law).
Given ${\mathbb{U}^B}=(\Omega^{\mathbb{U}^B}, (\F^{\mathbb{U}^B}_t), \mathbb{P}^{\mathbb{U}^B}, (W^{1,\mathbb{U}^B}_t), (W^{2,\mathbb{U}^B}_t), (B^{\mathbb{U}^B}_t),  (u^{\mathbb{U}^B}_t))  \in \mathbb{S} ^{1,2,B}$, we set $$ \widetilde{B}_t= \eta ^{-1} \int_0^t b(X^{\eta}_s, Q^{\e,\eta}_s,u ^{\mathbb{U}^B}_s)\, ds + {B}^{\mathbb{U}^B}_t \hbox{ and } \widetilde{W}^2= \e ^{-1/2}\int_0^t \rho(u^{\mathbb{U}^B}_s)\, ds + {W}_t^{2,\mathbb{U}^B}.$$
Hence $ (X ^{ \e, \eta,\mathbb{U}^B}, Q^{\e, \eta, \mathbb{U}^B })$ solves:
\begin{equation}
\begin{cases}
dX_t= AX_t\,dt  +  R(X_t) d{W}^{1,\mathbb{U}^B}_t + \eta \,  d \widetilde {B}_t, & X_0=x_0, \\ \\
 \e dQ_t= (BQ_t+ F(X_t,Q_t) dt +   \e^{1/2} G \, d\widetilde{W}^{2}_t, & Q^\e_0=q_0.
\end{cases}
\end{equation}
Now we associate to such forward system the backward equation 
\begin{equation}\label{BSDEsetting}
\begin{cases}
 -d Y_t= \psi^{\epsilon,\eta}(X^{\e, \eta,\mathbb{U}^B}_t, Q^{\e, \eta,\mathbb{U}^B}_t, {Z}^1_t, {Z}^2_t, \Xi_t ) dt-{Z}^1_t   d{W}^{1,\mathbb{U}^B}_t -  {Z}^2_t   d\widetilde{B}_t  - {\Xi}_t   d \widetilde{W}^{2}_t,  \\ \quad Y _1 \  = h(X_1).
\end{cases}
\end{equation}
Since the law of the solution does not depend on the particular setting, we get that $ Y_0 = Y^{\e, \eta}_0$, then rewriting \eqref{BSDEsetting} with respect to ${B}^{\mathbb{U}^B}$ and $ {W}^{2,\mathbb{U}^B}$, computing the expectation with respect to $\mathbb{P}^{\mathbb{U}^B}$ and recalling the definition of $\psi ^{\e, \eta}$, we get:
\begin{equation}
    Y _0= Y^{\e, \eta}_0 \leq    \E ^{{\mathbb{U}^B}}  \Big[\int_0^1 l (X^{\e, \eta, \mathbb{U}^B}_t, Q^{\e,\eta, \mathbb{U}^B}_t,u^{\mathbb{U}^B}_t)\, dt + h(X^{\eta, \mathbb{U}^B}_1)\Big]
 =J^{\e, \eta}(x_0,q_0,\mathbb{U}^B).
\end{equation}
Thus the proof is completed.
\finedim
\medskip

By simple considerations on the control problems 
we can prove the following uniform convergence.
\begin{theorem}{
Under \ref{A.1}---\ref{A.7}  we  have that:
\begin{equation}\label{limite epsilon fisso}
\lim_{\eta \to 0} \sup_{\e >0} |V^{\e,\eta}(x_0,q_0)- V^{\e, \eta}(x_0,q_0)|=0
\end{equation}}
\end{theorem}
{\bf Proof.}
First of all we recall, see  Lemma \ref{identita value function},  that $V^{\e}(x_0,q_0)= V^{\e,0}(x_0,q_0)$.

Moreover, by definition, $V^{\e}(x_0,q_0)=  \inf_{\mathbb{U}^B \in \mathbb{S} ^{1,2,B}} J^{\e,\eta}(x_0,q_0,\mathbb{U}^B)$, for every $\eta \geq 0$. 
Therefore, our claim follows if we prove that:
\begin{equation}
\lim_{\eta \to 0} \sup_{\e >0}  \sup _{\mathbb{U}^B \in \mathbb{S} ^{1,2,B}} |J^{\e,\eta}(x_0,q_0,\mathbb{U}^B)-J^{\e,0}(x_0,q_0,\mathbb{U}^B)|=0.
\end{equation}
We fix then $\e >0$ and $\mathbb{U}^B=(\Omega, (\F_t), \mathbb{P}, (W^1_t), (W^2_t), (B_t), (u_t)) \in \mathbb{S} ^{1,2,B}$ and consider:
$$\begin{cases}
d(X^{\e,\eta, \mathbb{U}^B}_t-X^{\e,0, \mathbb{U}^B}_t)= A(X^{\e,\eta,\mathbb{U}^B}_t-X^{\e,0,\mathbb{U}^B}_t) \, dt  +  [b(X^{\e,\eta,\mathbb{U}^B}_t,Q^{\e,\eta,\mathbb{U}^B}_t,u_t)-b(X^{\e,0,\mathbb{U}^B}_t,Q^{\e,0,\mathbb{U}^B}_t,u_t)]dt  & \\
\qquad\qquad\qquad\qquad + [R (X^{\e,\eta,\mathbb{U}^B}_t)-R(X^{\e,0, \mathbb{U}^B}_t)] dW^1_t + \eta\,  d B_t, & \\
X^{\e,\eta,\mathbb{U}^B}_0-X^{\e,0,\mathbb{U}^B}_0=0, \\ 
 d(Q^{\e,\eta,\mathbb{U}^B}_t-Q^{\e,0,\mathbb{U}^B}_t)= \varepsilon^{-1} B (Q^{\e,\eta,\mathbb{U}^B}_t-Q^{\e,0,\mathbb{U}^B}_t)dt+ 
 \e^{-1} [F(X^{\e,\eta,\mathbb{U}^B}_t,Q^{\e,\eta,\mathbb{U}^B}_t,u_t)-F(X^{\e,0,\mathbb{U}^B}_t,Q^{\e,0,\mathbb{U}^B}_t,u_t)]dt, & \\
Q^{\e,\eta,\mathbb{U}^B}_0-Q^{\e,0,\mathbb{U}^B}_0=0.&
\end{cases}$$
Taking into account the second equation we have by Hypothesis \ref{A.3} and standard estimates (see for instance \cite{DPZab2})
\begin{equation}\label{stimaQdaXintegrale}
    |Q^{\e,\eta,\mathbb{U}^B}_t-Q^{\e,0,\mathbb{U}^B}_t|\leq \frac{C}{\varepsilon} \int_0^t e^{-\displaystyle\frac{\mu}{\varepsilon}(t-s)} |X^{\e,\eta,\mathbb{U}^B}_s-X^{\e,0, \mathbb{U}^B}_s| \, ds,
\end{equation}
and consequently
\begin{equation}\label{stimaQdaXsup}
 \sup_{s\leq t}|Q^{\e,\eta,\mathbb{U}^B}_s-Q^{\e,0,\mathbb{U}^B}_s|\leq \frac{C}{\mu}
 \sup_{s\leq t}|X^{\e,\eta,\mathbb{U}^B}_s-X^{\e,0, \mathbb{U}^B}_s|
 \end{equation}
 where $C$ is a constant independent of $\epsilon$ and $\eta$ with value that can change from line to line.

As far as the first equation is concerned we have:
\begin{align*}
X^{\e, \eta, \mathbb{U}^B}_t  - X^{\e, 0,\mathbb{U}^B}_t =  &\int_0^t e^{(t-s) A}
[b (X^{\e, \eta, \mathbb{U}^B}_s, Q^{\e, \eta, \mathbb{U}^B}_s, u_s) - b (X^{\e, 0, \mathbb{U}^B}_s, Q^{\e, 0, \mathbb{U}^B}_s, u_s) ] \, ds  +  \int_0^t e^{(t-s) A}\eta  \, d B_s \\
&+\int_0^t e^{(t-s) A}
[R (X^{\e, \eta, \mathbb{U}^B}_s) - R (X^{\e, 0, \mathbb{U}^B}_s) ] \, d W^1_s.
\end{align*}
Thanks to standard estimates and the factorization method, see \cite{DpZfatt}, we get for $p> \frac{2}{1-2\gamma}$,  $\alpha \in (\frac{1}{p},\frac{1}{2} -\gamma)$, ($\gamma\in (0,1/2)$ is the constant appearing in Assumption \ref{A.1} and Assumption \ref{A.5}) and any $\rho \in [0,1]$:
\begin{align*}
& \E^{\mathbb{U}^B }\sup_{t \in [0,\rho]} |  X^{\e, \eta, \mathbb{U}^B}_t - X^{\e, 0,\mathbb{U}^B}_t |^p     \leq   C  \int_0^\rho [\E ^{\mathbb{U}^B}\!\!\!\sup_{s \in [0,r]} |X^{\e, \eta, \mathbb{U}^B}_s-  X^{\e, 0, \mathbb{U}^B}_s|^p + 
 \E^{\mathbb{U}^B}  |  Q^{\e, \eta, \mathbb{U}^B}_r - Q^{\e, 0,\mathbb{U}^B}_r |^p] \, dr   \\ 
&  + \int_0^\rho \E ^{\mathbb{U}^B}\Big[ \int_0^r (r-l)^{-2(\alpha+\gamma)} | X^{\e, \eta, \mathbb{U}^B}_l -  X^{\e, 0, \mathbb{U}^B}_l|^2\, dl \Big]^{p/2}\, dr +  |\eta|^p\\ & \leq  C  \left(1+ \left(\int_0^1 \sigma^{-2(\alpha+\gamma)}\, d\sigma \right)^{p/2} \right) \int_0^\rho \E ^{\mathbb{U}^B}\!\!\! \sup_{s \in [0,r]} |X^{\e, \eta, \mathbb{U}^B}_s-  X^{\e, 0, \mathbb{U}^B}_s|^p  \, dr + 
 \int_0^\rho \E ^{\mathbb{U}^B}  |  Q^{\e, \eta, \mathbb{U}^B}_r - Q^{\e, 0, \mathbb{U}^B}_r |^p \, dr  
   +   |\eta|^p.
\end{align*}
 Recalling \eqref{stimaQdaXsup} we also get:
$$\E^{\mathbb{U}^B}  \sup_{t \in [0,\rho]} |  X^{\e, \eta, {\mathbb{U}^B}}_t - X^{\e, 0, {\mathbb{U}^B}}_t |^p     \leq C \Big [\int_0^\rho \E ^{\mathbb{U}^B} \sup_{s \in [0,r]} |X^{\e, \eta, {\mathbb{U}^B}}_s-  X^{\e, 0,{\mathbb{U}^B}}_s|^p  \, dr  + |\eta|^p\Big ]; $$
and  applying the Gromwall Lemma
to $v(r)=: \E ^{\mathbb{U}^B}  \sup_{s \in [0,r]} |  X^{\e, \eta,  \mathbb{U}^B}_s - X^{\e, 0, \mathbb{U}^B}_s |^p $, we conclude
\begin{align*}
 &\E ^{\mathbb{U}^B}   \sup_{t \in [0,\rho]} |  X^{\e, \eta, \mathbb{U}^B }_t - X^{\e,0, \mathbb{U}^B }_t |^p   \leq  C |\eta|^p, \qquad \qquad \forall \e >0 
\end{align*}
and applying once again \eqref{stimaQdaXsup}
\begin{align*}
 &\E ^{\mathbb{U}^B} \sup_{t \in [0,\rho]} |  Q^{\e, \eta, \mathbb{U}^B}_t - Q^{\e,0, \mathbb{U}^B}_t |^p   \leq  C |\eta|^p, \qquad \qquad \forall \e >0 
\end{align*}
Finally if we consider the difference between the value functions:
\begin{align}
   | V^{\e, \eta} (x_0,q_0) - V^{\e,0} (x_0,q_0)| & \leq \sup _{  \mathbb{U}^B \in \mathbb{S}^{1,2,B} }\Big [\E ^{\mathbb{U}^B}  \int_0^1 |l(X^{\e, \eta, \mathbb{U}^B}_t,Q^{\e, \eta, \mathbb{U}^B}_t,u_t)- 
 l(X^{\e, 0,\mathbb{U}^B}_t,Q^{\e,0, \mathbb{U}^B}_t,u_t)| \, dt \nonumber\\ & \qquad\qquad + \E ^{\mathbb{U}^B} |h(X^{\e, \eta, \mathbb{U}^B}_1)- h(X^{\e,0, \mathbb{U}^B}_1)| \Big ]\\ &  \leq  C  \sup _{  \mathbb{U}^B \in \mathbb{S}^{1,2,B} }
  \E^{\mathbb{U}^B}   \Big[  \sup_{t \in [0,1]}|  X^{\e, \eta, \mathbb{U}^B }_t - X^{\e,0, \mathbb{U}^B }_t | +  \int_0^1 |  Q^{\e, \eta, \mathbb{U}^B }_t - Q^{\e, 0, \mathbb{U}^B }_t |  \, dt \Big] \nonumber \\& \leq C |\eta|, \qquad \forall \e >0.\nonumber
\end{align}
Thus our claim holds.
\finedim

\begin{theorem}{\label{ teo_eta fisso}
For every fixed $x\in H$ and $z\in H^*$, let us consider the following ergodic control problem with {\em state equation} in $K$, driven by an arbitrary $\Xi$-valued cylindrical Wiener process $(\hat{{W}})$, with control $\beta: [0,\infty[ \times \Omega \rightarrow U$ varying in the set $\mathcal{U}_{\infty}$ of progressively measurable $U$-valued processes with respect to the natural filtration of  $(\hat{{W}})$. 

\begin{equation}\label{ergfrozenstate}
d \hat{Q}^{\beta}_s=B\hat{Q}^{\beta}_sds+ F(x, \hat{Q}^{\beta}_s)\, d s + G\rho(\beta_s) d s +  G d\hat {W}^2_s, \qquad  \hat{Q}^{\beta}_0=0
\end{equation}
and {\em ergodic cost functional}:

\begin{equation}\label{ergfrozencostT}
 \check{J}(x,z,\beta)=\liminf_{\delta \to 0} \hat{\E} \, \delta \int_0^{\frac{1}{\delta}}  [z \, b(x, \hat{Q}^{\beta}_s,\beta_s)+ l(x,\hat{Q}^{\beta}_s,\beta_s)] d s.
\end{equation}

Let $\lambda(x,z)$ be is the value function of the above ergodic control problem, that is:
\begin{equation}
\label{lambda}
\lambda(x,z)=\inf_{\beta\in \mathcal{U}^{2}}  \check{J}(x,z,\beta).
\end{equation}
Under \ref{A.1}---\ref{A.7}  we  have that $\lambda(x,\,\cdot\,)$ is concave moreover:
\begin{equation}\label{proprietàdilambda}
\begin{array}{l}
     | \lambda(x,z)-\lambda(x,z')|\leq M |z-z'| \\
      | \lambda(x,z)-\lambda(x',z)|\leq L (1+|z|) |x-x'| \\
      |\lambda(x,z)|\leq M (1+|z|)
\end{array}
   \end{equation}

Moreover for every $\eta>0$:
\begin{equation}\label{limite eta fisso}
\lim_{\e \to 0}  Y^{\e,\eta}_0= \lim_{\e \to 0} V^{\e,\eta}(x_0,q_0)= Y^\eta_0
\end{equation}
where $Y^\eta_0$ is defined as part of the solution to the \textit{reduced} BSDE (for the definition of $X^{\eta }$ see \eqref{lentoveloceottimo}).
\begin{equation}\label{LimitEquation}
\begin{cases}
-d Y_t  = \lambda({X}^\eta_t,\eta^{-1}{Z}^{2}_t)\,dt - \, {Z}^{1}_t\, d\bar{W}^1_t  - \, {Z}^{2}_t\, d\bar{B}_t,   \\
 \quad {Y}_1=h(X ^\eta_1).
\end{cases} 
\end{equation}
where $Y^\eta \in L^2_{\bar{\mathcal{F}}^{1,B}}(\Omega;C([0,1];\R))$, $Z^{1,\eta} \in L^2_{\bar{\mathcal{F}}^{1,B}}(\Omega\times [0,1]; \Xi^* )$ and $Z^{2,\eta} \in L^2_{\bar{\mathcal{F}}^{1,B}}(\Omega\times [0,1]; H^* )$.}
\end{theorem}
{\bf Proof:}
The proof follows the one in \cite[Theorem 5.4]{GuaTess2019AMO}.
The only point to check is the discretization procedure of the forward component $X^\eta$ since now in the limit equation \eqref{LimitEquation} a multiplicative noise appears.
As in  \cite[Theorem 5.4]{GuaTess2019AMO} we introduce, for all $N \in \mathbb{N}$ a partition $\pi$ of $[0,1]$, $\pi= \{ \frac{k}{2^N}: k: 0, \dots, 2^N-1\}$ and we define a sequence $ X^{\eta, N} $
\begin{equation}
    X^{\eta,N}_t = \sum_{k=0}^{2^N-1} X^\eta_{\frac{k}{2^N}}I_{[\frac{k}{2^N}, \frac{k+1}{2^N}[}  (t) +   X^\eta_{1} \delta_1(t), \qquad t \in [0,1]
\end{equation}

We need to arrive to the following
\begin{equation}
  \lim_{N \to \infty}  \bar{\E} \sup_{t\in [0,1]} |X^{\eta, N}_t - X^\eta_t |^4 =0,
    \end{equation}
in order to exploit the procedure of \cite[Theorem 5.4]{GuaTess2019AMO}.

First of all, using the factorization method,  see  \cite{FuTessAOP2002} and \cite{DpZfatt}, we can find, for any $p> \frac{2}{1-2\gamma}$ (and consequently for any $p\geq 1$),   a constant $C_p$ independent of $\eta$ such that:
\begin{equation}\label{stimaXeta}
    \bar{\E} \sup_{t\in [0,1]} | X^\eta_t |^p  \leq C_p (|x_0|^p   + |\eta|^p).
\end{equation}

For $t \in [\frac{k}{2^N}, \frac{k+1}{2^N}[$ we evaluate the difference:
\begin{align}\label{diffatfixed}
    X^{\eta,N}_t -  X^{\eta}_t  = \int_{\frac{k}{2^N}}^t e^{(t-s)A }R (X^\eta_s) \, d \bar{W}^1_s+ \int_{\frac{k}{2^N}}^t  e^{(t-s)A}\eta \, d \bar{B}_s.
\end{align}
again using the factorization method \cite{DpZfatt} for any $p> \frac{2}{1-2\gamma}$ any $\alpha \in (\frac{1}{p}, \frac{1}{2}-\gamma)$ we can find a constant $M$, that depends on $\alpha$ but not on $N$ such that, for $ q= \frac{p}{p-1}$, we get:
\begin{align}\label{stimadiffNeta}
     &\bar{\E} \sup_{t\in [\frac{k}{2^N},\frac{k+1}{2^N}]} |X^{N,\eta}_t - X^\eta_t |^p \leq \nonumber \\&  M \Big(\int_{0}^{\frac{1}{2^N}} s^{{q}(\alpha-1)}\, dr \Big) ^{p/q}   \Big(\int_{0}^{1} (s-r)^{-2(\alpha+\gamma)}\, dr \Big)^{p/2} \int_0^1(1+ \bar{\E} \sup_{t \in [0,1]} |X^\eta_t|^p + |\eta| ^p) \, dt, \qquad \forall k\in \{0,\cdots , 2^N -1\}
    \end{align}
    
    For the reader convenience we write some details regarding the first term at the R.H.S. in  \eqref{diffatfixed}.
    \begin{align*}
 \int_{\frac{k}{2^N}}^t e^{(t-s)A }R (X^\eta_s) \, d \bar{W}^1_s &= \frac{1}{B(\alpha,1-\alpha)} \int_{\frac{k}{2^N}}^t e^{(t-r)A } (t-r)^{1-\alpha} \int_\frac{k}{2^N}^r e^{(r-s)A } (r-s)^{-\alpha}   R (X^\eta_s) \, d \bar{W}^1_s\, dr \\ &=  \frac{1}{B(\alpha,1-\alpha)} \int_{\frac{k}{2^N}}^t e^{(t-r)A } (t-r)^{1-\alpha} Y(r)\, dr 
 \end{align*}
    where by $B(\alpha,1-\alpha)$ we denote the normalization constant of the beta distribution and by
    $Y(r)$ the random variable $Y(r):=\int_\frac{k}{2^N}^r e^{(r-s)A } (r-s)^{-\alpha}   R (X^\eta_s) \, d W^1_s$.

 Thus for any $p > \frac{2}{1-2\gamma}$ and $  \frac{1}{2}-\gamma > \alpha > \frac{1}{p}$    
 \begin{align*}
 &\bar{\E} \sup_{t \in [\frac{k}{2^N},\frac{k+1}{2^N}]} \Big| \int_{\frac{k}{2^N}}^t e^{(t-s)A }R (X^\eta_s) \, d W^1_s \Big| ^ p\leq   {\Big (\frac{M_A e^{\omega_A}}{B (\alpha,1-\alpha)} \Big) ^p }\Big (\int_{\frac{k}{2^N}}^{\frac{k+1}{2^N}} (t-r)^{(1-\alpha)q} \, dr \Big) ^{p/q} \bar{ \E }  \int_{\frac{k}{2^N}}^{\frac{k+1}{2^N}} |Y(r)|^ p\, dr    \\ & \leq  {C} \Big(\int_{0}^{\frac{1}{2^N}} s^{{q}(\alpha-1)}\, dr \Big) ^{p/q} \int_{\frac{k}{2^N}}^{\frac{k+1}{2^N}} \bar{ \E  }\left| \int_\frac{k}{2^N}^r e^{(r-s)A } (r-s)^{-\alpha}   R (X^\eta _s ) \, d W^1_s\right |^ p \, dr  \\&\leq{C}  \Big(\int_{0}^{\frac{1}{2^N}} s^{{q}(\alpha-1)}\, dr \Big) ^{p/q}     \int_{\frac{k}{2^N}}^{\frac{k+1}{2^N}} \Big(\int_\frac{k}{2^N}^r (r-s)^ {-2\alpha-2\gamma}\bar{ \E }| R(X^\eta _s)|^ 2_{L(H)} \, ds \Big)^ {p/2} \, dr \\
 & \leq C  \Big(\int_{0}^{\frac{1}{2^N}} s^{{q}(\alpha-1)}\, dr \Big) ^{p/q} \Big(\int_0^1 \sigma^ {-2\alpha-2\gamma} \, d\sigma \Big)^ {p/2}   \Big(  \int_{0}^{1}  \bar{\E} \sup_{r \in [0,1]} (1+|X^\eta_r|^p)   \, dr \Big)
 \end{align*}   
Where ${C}$ may change from line to line but is always independent of $N$. Thus thanks to \eqref{stimaXeta} and \eqref{stimadiffNeta} we get the thesis.
\finedim.

\section{Interchanging limits}

We now prove general result allowing us to interchange the limit with respect to $\e$ and the one with respect to $\eta$ as well.

\begin{Theorem}\label{lemma_scambio_limiti}
Let $v^{\eta}(x_0)=Y^{\eta}_0$ (see \eqref{LimitEquation} for the definition of $Y^{\eta}$) then, for all $x_0\in H$ and $q_0\in K$, it holds:
$$\lim_{\e \rightarrow 0}V^{\e}(x_0, q_0)=\lim_{\eta\rightarrow 0} v^{\eta}(x_0):= V(x_0).$$
Moreover $V^{\e}$ is Lipschitz uniformly with respect to $\e$, $v^{\eta}$ is Lipschitz uniformly with respect to $\eta$ and $V$ is Lipschitz.
\end{Theorem}
{\bf Proof. } Since, fixed $x_0\in H$ the sequence {$v^{\eta}(x_0)$} is bounded (see \eqref{limite eta fisso} and \eqref{cost_epsilon_eta}) then there exists a sequence $\eta_n \searrow 0$ (depending on $x_0$ but we omit this information in the notation since it is not relevant here)  such that the sequence $v^{\eta_n}(x_0)$ converges to a limit that we denote by $V(x_0)$. By standard adding and subtracting
\begin{equation}\label{addingandsubtracting}
|V^{\e,0}(x_0,q_0)-V(x_0)|\leq | V^{\e,0}(x_0,q_0) -V^{\e,\eta_n}(x_0,q_0)|+
|V^{\e,\eta_n}(x_0,q_0)-v^{\eta_n}(x_0)|+
|v^{\eta_n}(x_0)-V(x_0)|
    \end{equation}
Fix $\delta >0$. By  \eqref{limite epsilon fisso} there exists $n_{\delta}$ such that 
$$|V^{\e,\eta_n}(x_0, q_0)-V^{\e,0}(x_0,q_0)|+
|v^{\eta_n}(x_0)-V(x_0)|\leq \delta\quad \forall \e>0,\; \forall n\geq n_{\delta}$$
We fix an arbitrary $\bar n \geq  n_{\delta}$ and notice that by Theorem \ref{ teo_eta fisso} there exists $\e_{\delta}>0$ such that 
$$|V^{\e,\eta_{\bar n}}(x_0.q_0)-v^{\eta_{\bar n}}(x_0)|<\delta, \qquad \forall \e\in (0,\e_{\delta}).$$
It is then straightforward to conclude:
$$\lim_{\e\rightarrow 0} V^{\e}(x_0.q_0)=V(x_0)$$
Moreover, by the same argument,  from any sequence $\hat{\eta}_n \searrow 0 $ we can extract a subsequence $\hat{\eta}_{n_k} \searrow 0 $ such that
 $$\lim_{k \rightarrow \infty } v^{\hat{\eta}_{n_k}}(x_0) =  \lim_{\e \rightarrow 0}V^{\e}(x_0, q_0)= V(x_0)$$
 and this implies that
 $$\lim_{\eta \rightarrow 0} v^{\eta}(x_0) =  \lim_{\e \rightarrow 0}V^{\e}(x_0, q_0)= V(x_0).$$
 It lasts to show that $V$ is Lipschitz continuous. Clearly it is enough to show that $V^{\e,\eta}$ is Lipschitz (with respect to $x_0$) uniformly in $\e >0$ and $\eta >0$. 
 
 First we notice that, by the definition of $V^{\e,\eta}(x,q) $ (we here indicate dependence of the solution of equation \ref{lentovelocecontrollata} on initial data $(x,q)$ and $(x',q')$) and hypothesis \eqref{A.7} :
\begin{align*}
    & |V^{\e,\eta}(x,q) - V^{\e,\eta}(x',q')|\leq \sup_{\mathbb{U} ^B \in \mathbb{S} ^{1,2,B}} \Big| \E ^{\mathbb{U} ^B} \int_0^1 l(X^{\e,\eta, \mathbb{U} ^B , x, q}_s, Q^{\e,\eta, \mathbb{U} ^B , x,q}_s,u_s)-l(X^{\e,\eta, \mathbb{U} ^B , x', q'}_s, Q^{\e,\eta, \mathbb{U} ^B , x',q'}_s,u_s)\, ds \Big|\\
    &  + | \E^{\mathbb{U} ^B}  [h(X^{\e,\eta, \mathbb{U} ^B , x, q}_1)-  h (X^{\e,\eta, \mathbb{U} ^B , x', q'}_1)]|, \\
    & \leq L [ \sup_{\mathbb{U} ^B \in \mathbb{S} ^{1,2,B}}
    \E^{\mathbb{U} ^B}   \sup_{s\in [0,1]}   |X^{\e,\eta, \mathbb{U} ^B, x, q}_s- X^{\e,\eta, \mathbb{U} ^B, x', q'}_s| + \E ^{\mathbb{U}^B}  
    \int_0^1 |Q^{\e,\eta, \mathbb{U}^B, x, q}_s- Q^{\e,\eta, \mathbb{U}^B, x', q'}_s| \, ds].
\end{align*}
 Then, arguing as in Theorem \ref{limite epsilon fisso}, we can use the factorization method and the dissipativity condition to get that  $\forall r \in [0,1]$
 \[ \E ^{\mathbb{U} ^B}\sup_{s\in [0,r]}   |Q^{\e,\eta, {\mathbb{U} ^B}, x, q}_s- Q^{\e,\eta, {\mathbb{U} ^B}, x', q'}_s| \leq C (  \E ^{\mathbb{U} ^B} \sup_{s\in [0,r]}   |X^{\e,\eta, {\mathbb{U} ^B}, x, q}_s- X^{\e,\eta, {\mathbb{U} ^B}, x', q'}_s|   + |q-q'|), \]
 and consequently that
\[ \E ^{\mathbb{U} ^B} \sup_{s\in [0,1]}   |X^{\e,\eta, {\mathbb{U} ^B}, x, q}_s- X^{\e,\eta, {\mathbb{U} ^B}, x', q'}_s| \leq C ( |x-x'| + |q-q'|), \]
 for some positive constants $C$ independent from $u$   $\e$ and $\eta$.
 \finedim{}

\section{Main characterizations}



For $\eta$ fixed we consider the limit system (starting now at an arbitrary time $s \in [0,1]$ from state $x\in H$ and written with respect to the reference setting).
\begin{equation}\label{forbaceta}
\begin{cases}
d X_t =   AX_t \, dt + R(X_t)\, d \bar{W}^1 _t + \eta \, d \bar{B}_t, \qquad t \in [s,1], \\
-d Y _t  = \lambda({X}_t,\eta^{-1}{Z}^2_t)\,dt - \, {Z}^1_t\, d \bar{W}^1_t + \, {Z}^2_t\, d \bar{B}_t,   \\
X_s= x, \quad {Y} _1=h(X_1).
\end{cases} 
\end{equation}
Again by standard BSDE theory equation \eqref{forbaceta} has a unique solution
 $(Y^{ \eta,s,x}, Z^{1,\eta,s,x}, Z^{2, \eta,s,x}, \Xi^{ \eta,s,x})$ with   $Y^{ \eta,s,x}\! \in L^2_{\bar{\mathcal{F}}^{1,B}}(\Omega;C([s,1];\R))$, $Z^{ 1,\eta,s,x} \!\in L^2_{\bar{\mathcal{F}}^{1,B}} (\Omega\times [s,1];\Xi^*)$, $Z^{ 2,\eta,s,x} \!\in L^2_{\bar{\mathcal{F}}^{1,B}} (\Omega\times [s,1];H^*)$ and  $\Xi^{ \eta,s,x} \!\in L^2_{\bar{\mathcal{F}}^{1,B}}(\Omega\times [s,1];\Xi^*)$,
where $\lambda(x,z)$ is defined in \eqref{lambda} (notice that in this section we need to indicate the dependence on the initial time $s$ and state $x$).

The above system follows within the framework of  \cite[Theorem 4.1]{BismutElworthy} (indeed  assumptions  2.1, 3.1-3.7 and 4.1  
are verified and hypothesis 3.4 can be relaxed as  pointed out in \cite[pag. 443]{BismutElworthy}).
Thus if we denote with $v^\eta (s,x)= Y^{\eta, s,x}_s $, we have that 
$$ (Z^{1,\eta, s,x}_t,Z^{2,\eta, s,x}_t)=\left(\nabla_x v^{\eta}(t,X^{\eta, s,x}_t )R (X^{\eta, s,x}_t), \eta \nabla_x v^{\eta}(t,X^{\eta, s,x}_t ) \right) $$
In particular $
 Z^{2,\eta,s,x}_t = \eta \nabla_x v^\eta(t, X^{\eta, s,x}_t )$,  $ \mathbb{P} \times ds$ -almost surely.

Taking into account the representation in Lemma \ref{lemma 4.3} (with initial time $s$ instead of 0) and proceeding as in the proof of Theorem \ref{lemma_scambio_limiti}, we get that
$v^\eta(s, \cdot)$ is Lipschitz uniformly with respect to $s\in[0,1]$ and $\eta >0$. Thus (see also 
 \cite[Theorem 4.1]{BismutElworthy} ) we have that:
\begin{equation}\label{stimauniffacile}
  \displaystyle\Big|{ Z^{2,\eta,s,x}_t}\Big|  _{H^*} \leq a \eta , \qquad d\mathbb{P} \times dt -a.e..
\end{equation}
where $a\in \mathbb{R}^+$ is independent on $\eta$, $s$ and $x$.

We now introduce the following function, $\tilde \lambda$ defined, for a constant $k>M$  large enough, by:
\begin{equation}\label{tildelambda}
    \tilde \lambda (x,z): = \lambda (x,z) \wedge [-(M+1) |z| + \kappa]
\end{equation}
By \eqref{proprietàdilambda} we get
\begin{equation*}
    \tilde \lambda (x,z): =  \begin{cases}
    \lambda (x,z), &  \hbox{ if } |z| \leq \displaystyle \frac{\kappa-M}{2M+1}, \\
    \kappa-(M+1) |z|, & \hbox{ if }|z| \geq \displaystyle k+M.
    \end{cases}
\end{equation*}
Choosing $k$ large enough we can  assume that $(k-M)/(2M+1)>a$ so that $\tilde \lambda (x,z)=\lambda (x,z)$ when $|z|\leq a$.  Moreover   $\tilde{\lambda}$ remains concave being the minimum of concave functions and, by  \eqref{proprietàdilambda}, is Lipschitz with respect  to $z$, uniformly with respect to $x$ with Lipschitz constant $M+1$.

Finally we have:
\begin{equation*}
   | \tilde \lambda (x,z) -  \tilde \lambda (x',z)|\; \begin{cases}=
     | \lambda (x,z) -   \lambda (x',z)|, &  |z| \leq a, \\ 
     \leq   | \lambda (x,z) -   \lambda (x',z)|, & a <|z| < \kappa+M \\
    =0, & |z|
    \geq\kappa+M.
    \end{cases}
\end{equation*}
therefore, thanks again to \eqref{proprietàdilambda}, $\tilde \lambda$  is Lipschitz  continuous with respect to $x$, uniformly w.r.t. $z$, with  Lipschitz  constant equal to $\tilde{L}=L(1+\kappa+M)$.

Taking into account \eqref{stimauniffacile}, system \eqref{forbaceta}  can be written, replacing $\lambda$ by $\tilde\lambda$ and choosing $s=0$, as follows:
\begin{equation}\label{forbacetatilde}
\begin{cases}
d X_t =   AX_t \, dt + R (X_t)\, d \bar{W}^1 _t + \eta \, d \bar{B} _t , \qquad t \in [0,1], \\
-d Y _t  = \tilde{\lambda}({X}_t,\eta^{-1}{Z}^{2}_t)\,dt - \, {Z}^{1}_t\, d\bar{W}^1_t  - \, {Z}^{2}_t\, d\bar{B}_t,   \\
X_0= x, \quad {Y} _1=h(X_1).
\end{cases} 
\end{equation}

We denote by  $\tilde\lambda_*$ the Legendre transform  of $\tilde\lambda$, that is, for $x$ and  $\alpha$ in $H $ (recall that $\tilde{\lambda}$ is concave, this justifies the negative signs):
\begin{equation}\label{tildelambdastar}
    \tilde \lambda _*(x,\alpha): =  \inf_{z \in H^*} \{-z\alpha -\tilde \lambda {(x,z)} \}
\end{equation}
It turns out that     $\tilde \lambda _*$ is Lipschitz continuous with respect to $x$, uniformly w.r.t. $\alpha$, as well. Indeed:
\begin{equation}\label{LegendreLip}
|\tilde \lambda _*(x,\alpha) - \tilde \lambda _*(x',\alpha)| \leq \sup _{z \in H^*}
|\tilde \lambda(x,z) - \tilde \lambda (x',z)| \leq \tilde{L} \,  |x-x'|, \qquad \forall x ,x', \alpha\in H.
\end{equation}

Moreover taking into account Lipschitzianity with respect to $z$ of $\tilde{\lambda}$ we get:
$$ \tilde \lambda_*(x,\alpha)= -\infty \; \hbox{ if } |\alpha|>M+1 $$

That yields the following simplification in the Fenchel duality:
\begin{equation}\label{dominiolimitato}
    \tilde \lambda(x,z): =  \inf_{\alpha \in H } \{-z\alpha -\tilde \lambda_* {(x,\alpha)} \}
    =  \inf_{\alpha \in H: |\alpha|\leq M+1} \{-z\alpha -\tilde \lambda_* {(x,\alpha)} \} 
\end{equation}
The solution $(Y^\eta)$ can then be represented by a reduced control problem that has the needed regularity to eventually allow the passage to the limit as $\eta \rightarrow 0$ giving the final representation of $\lim_{\e \rightarrow 0} V^{\e}$ .

We denote by $\mathcal{U}^{1,B}_{H}$ the set of all processes $(\alpha_t)_{t\in [0,1]}$  taking values in the ball $\{\alpha\in H: |\alpha|\leq M+1\}$ and being progressively measurable with respect to the  filtration $(\bar{\mathcal{F}}^{1,B})$.

\begin{lemma}\label{caratt-Y-eta}


We have:
\begin{equation}\label{rappreYeta}
    Y^{\eta} _t =  \inf_{\alpha  \in \,\mathcal{U}^{1,B}_{H}} \E^{\alpha}  \biggl( h(X^{\eta}_{1}) - \int_{t}^1  \tilde \lambda_*(X^{\eta}_{\ell},{\alpha}_{\ell}) d\ell \bigg | \mathcal{F}^{1,B}_t \biggr),
\end{equation}
where 
${\E}^{{\alpha}}$ denotes the mean value with respect the probability $\mathbb{P}^{{\alpha}}$ under which
$$
\left ( \bar{W}^1_t,\displaystyle \int_0^t \frac{{\alpha}_{\ell}}{\eta} d\ell + \bar{B}_t \right) : = ( \bar{W}^1_t, \bar{B}^\alpha_t) $$ is a Wiener process.

Notice that, with respect to $(\bar{W}^1,\bar{B}^{{\alpha}})$ process $(X^{\eta})$ solves the controlled stochastic differential equation:
\begin{equation}\label{Xetap}
    dX_t=
AX_t dt  - \, {\alpha}_t dt  +  R (X_t) d\bar{W}^1_t + \eta \bar{B}^\alpha_t, \quad X_0=x.
\end{equation}

\end{lemma}

\textbf{Proof:} 
To start with we point out that in the \eqref{dominiolimitato} the infimum can be restricted  to a bounded subset of $H$, as a consequence the choice of controls $\alpha$ in \eqref{Xetap} can be restricted to bounded (by $M+1$) controls and we are allowed to apply Girsanov transform to see perturbation by $\alpha$ as a change of probability.

Taking into account \eqref{dominiolimitato}, equation \eqref{forbacetatilde} evaluated at its solution $(X^\eta, Y^\eta, Z^{1,\eta},Z^{2,\eta})$ yields:
\begin{equation} \label{rel_fond_Y}
\begin{array}{rcl}
{Y}^\eta_t &=& \displaystyle h(X^ \eta_1)+\int_t^1 \tilde\lambda(X^ \eta_\ell,\eta^{-1}{Z}_{\ell}^{2, \eta})ds -\int_t^1 {Z}^{1,\eta}_\ell d\bar{W}^1_\ell -\int_t^1 {Z}^{2,\eta}_\ell d\bar{B}_\ell \\ &\leq& h(X^\eta_1) \displaystyle -\int_t^1 \left(\eta^{-1}{Z}^{2, \eta}_\ell{\alpha}_{\ell}+\tilde\lambda_*(X^{\eta}_\ell, {\alpha}_\ell)\right) d\ell -\int_t^1 Z^{1,\eta} _\ell d\bar{W}^1_\ell -\int_t^1 {Z}^{2, \eta}_\ell d\bar{B}_\ell .\end{array}     
\end{equation}

and by the definition of $(\bar{B}^{\alpha})$:
$${Y}^{\eta}_t\leq  h(X^{\eta} _1) \displaystyle -\int_t^1 \tilde\lambda_*(X^{\eta}_\ell, {\alpha}_\ell) d\ell -\int_t^1 Z^{1,\eta}_\ell d\bar{W}^1_\ell -\int_t^1 {Z}^{2,\eta}_\ell d\bar{B}^\alpha_\ell,$$
which shows that for all $\alpha\in \mathcal{U}_H^{1,B}$, 
\begin{equation*}
{Y}^{\eta}_{t} \leq \bar{\E}^{\alpha} \biggl( h(X^\eta_{1}) - \int_{t}^1 \tilde \lambda_*(X^{\eta}_{\ell},{\alpha}_{\ell}) d \ell \bigg | \bar{\mathcal{F}}^{1,B}_t\biggr).
\end{equation*}

Conversely, by measurable selection, we may choose a minimizing sequence of controls, 
$(\bar{\alpha}^n)_{n\in \mathbb{N}} \subset \mathcal{U}^{1,B}$,  such that, for all $\ell\in [0,1]$,  $\mathbb{P}$-a.s.:
 \begin{equation}\label{minimizingp}
  - \frac{Z^{2,\eta}_{\ell}}{\eta}\bar{\alpha}_{\ell }^n - \tilde\lambda_*(X^{\eta}_{\ell}, \bar{\alpha}_{\ell}^n)-1/n\leq \tilde\lambda \left(X^{\eta}_{\ell},\frac{Z^{2,\eta}_{\ell}}{\eta}\right).
 \end{equation}
Proceeding as in \eqref{rel_fond_Y} and taking into account \eqref{minimizingp} to obtain the reverse inequality we get:
$${Y}^{\eta}_t\geq  h(X_1^{\eta})  -\int_t^1 \left(\frac{{Z}^{2,\eta}_\ell}{\eta} \bar{\alpha}^n_\ell+\tilde \lambda_*(X^{\eta}_\ell, \bar{\alpha}^n_{\ell})+\frac{1}{n}\right) d \ell -\int_t^1 {Z}^{1,\eta}_\ell d\bar{W}^1_s- \int_t^1 {Z}^{2,\eta}_\ell d\bar{B}_\ell,
$$ 
and rewriting the above in terms of $\bar{B}^{{\bar{\alpha}}^n}$:
\begin{equation*}
\begin{split}
{Y}^{\eta}_{t} +\frac{1-t}{n} &\geq h(X^\eta_1)- \int_t^1  \tilde\lambda_*(X^{\eta}_{\ell},\bar{\alpha}^n_{\ell}) d \ell  - \int_t^1 {Z}^{1,\eta}_\ell d\bar{W}^1_\ell- \int_t^1 {Z}^{2,\eta}_\ell d\bar{B}^{\bar{\alpha}^n}_\ell.
\end{split}
\end{equation*}
Therefore we can conclude that:
$$
{Y}^{\eta}_{t} + 1/n  \geq \bar{\E}^{{\bar{\alpha}^n}} \biggl( h(X^{\eta}_{1}) - \int_{t}^1 \tilde{\lambda}_*(X^{\eta}_{\ell},\bar{\alpha}^n_{\ell}) d\ell \bigg | \bar{\mathcal{F}}^{1,B}_t\biggr) 
$$
and the claim is proved. \finedim

To complete the circle and give sense to the limit as $\eta \rightarrow 0$ we just have to come back to a control representation of $Y^{\eta}_0$ that does not relay on Girsanov transform (meaningless if $\eta=0$). This was already done for a different control problem but by general techniques in Lemma \ref{lemma 4.3}.
\medskip

Namely  if we denote  by $\mathbb{S} ^{1,B}_H$ the class of $7$-uples  $\mathbb{U}^B_H=(\Omega,\F, (\F_t), \mathbb{P}, (W^1_t),  (B_t), (\alpha_t))$, similar to the settings 
$\mathbb{S} ^{1,2,B}$ defined in Section \ref{sec-statement} with the only differences are that here the noise $(W^2)$ is omitted and $(\alpha)$ is now any $(\F_t)$-progressive process with values in the closed ball $\{x\in H : |x|\leq M+1\}$ (the subscript $H$ in the notation $\mathbb{S}_H$ indicates that we are now considering $H$-valued controls). It holds:
\begin{lemma}
\begin{equation}\label{cataratt Y eta}
    Y^{\eta}_0 = \inf_{  \mathbb{U}^{B}_H  \in \mathbb{S} ^{1,B}_H} 
    \bar{\E}^{  \mathbb{U}^{B}_H } \biggl( h(X^{\eta,  \mathbb{U}^{B}_H }_{1}) - \int_{t}^1  \tilde \lambda_*(X^{\eta,  \mathbb{U}^{B}_H }_{\ell},{\alpha}_{\ell}) d \ell  \biggr)
\end{equation}
where given  $   \mathbb{U}^{B}_H\in \mathbb{S}^{1,B}_H$ as above 
$X^{\eta,   \mathbb{U}^{B}_H} $ solves:

\begin{equation}\label{eq_X_eta}
    d{X}_s= A{X}_s ds  - \, {\alpha}_s ds  +  R ({X}_s) dW^1_s + \eta \, d B_t, \quad {X}_0=x_0.
\end{equation}
\end{lemma}
\Dim

The proof follows exactly as in the cited Lemma \ref{lemma 4.3}. \finedim

We are now able to prove the main result of the paper namely the characterization of the limit, as $\e\rightarrow 0$, of the value function $V^{\e}(x_0,q_0)$ of the original control problem in terms the value function of a \textit{reduced} control problem on a \textit{reduced} state space.

Let $\mathbb{S} ^{1}_H$ the class of $6$-uples  $\mathbb{U}_H=(\Omega,\F, (\F_t), \mathbb{P}, (W^1_t), (\alpha_t))$ identical to the ones in $\mathbb{S} ^{1,B}_H$ with the only difference that $(B)$ is not present.

\begin{theorem}\label{main-caratt}
It holds:
\begin{align*}
\lim_{\e\rightarrow 0} V^{\e}(x_0,q_0) ={V}(x_0) &=\inf_{ \mathbb{U}_H \in \, \mathbb{S}^{1}_H} \mathbb{E}^{\mathbb{U}_H}   \left( h({X}^{\mathbb{U}_H}_{1}) - \int_{0}^1  \tilde \lambda_*({X}^{\mathbb{U}_H}_s,{\alpha}_{s}^{\mathbb{U}_H}) ds \right).
    \end{align*} 
     where, given  $\mathbb{U}_H\in \mathbb{S}^{1}_H$ as above,
$X^{\eta,   \mathbb{U}_H} $ solves the state equation: 
$$d{X}_s= A{X}_s ds  - \, {\alpha}_s ds  +  R ({X}_s) dW^1_s, \quad {X}_0=x_0.$$
\end{theorem}

\textbf{Proof:} 
First we notice that as in Lemma \ref{identita value function} it holds:
$$\inf_{ \mathbb{U}_H \in \, \mathbb{S}^{1}_H} \mathbb{E}^{\mathbb{U}_H}   \left( h({X}^{\mathbb{U}_H}_{1}) - \int_{0}^1  \tilde \lambda_*({X}^{\mathbb{U}_H}_s,{\alpha}_{s}^{\mathbb{U}_H}) ds \right)=\inf_{ \mathbb{U}^{B}_H \in \, \mathbb{S}^{1,B}_H} \mathbb{E}^{\mathbb{U}^{B}_H }   \left( h({X}^{0,\mathbb{U}^{B}_H }_{1}) - \int_{0}^1  \tilde \lambda_*({X}^{0,\mathbb{U}^{B}_H }_s,{\alpha}_{s}^{\mathbb{U}^{B}_H }) ds \right).$$
where, we recall ${X}^{0,\mathbb{U}^{B}_H }$ solves equation \eqref{eq_X_eta} with $\eta=0$.

Thus, by \eqref{cataratt Y eta} and Th. \ref{lemma_scambio_limiti}  it is enough to prove that, if $\eta \rightarrow 0$:
$$ \inf_{  \mathbb{U}^{B}_H  \in \mathbb{S} ^{1,B}_H} 
    \bar{\E}^{  \mathbb{U}^{B}_H } \biggl( h(X^{\eta,  \mathbb{U}^{B}_H }_{1}) - \int_{t}^1  \tilde \lambda_*(X^{\eta, \mathbb{U}^{B}_H }_{\ell},{\alpha}_{\ell}) d \ell  \biggr)\rightarrow  \inf_{  \mathbb{U}^{B}_H  \in \mathbb{S} ^{1,B}_H} 
    \bar{\E}^{  \mathbb{U}^{B}_H } \biggl( h(X^{0,  \mathbb{U}^{B}_H }_{1}) - \int_{t}^1  \tilde \lambda_*(X^{0,  \mathbb{U}^{B}_H }_{\ell},{\alpha}_{\ell}) d \ell  \biggr)$$

 Fix any $   \mathbb{U}^{B}_H = (\Omega, (\F_t), \mathbb{P}, (W^1_t), (B_t),  (\alpha_t))  \in \mathbb{S} ^{1,B}_H$.
Thanks to the Lipschitzianity of $h$ (see hypotheses \ref{A.7}) and of $\tilde{\lambda}_*$  we easily have that for a suitable constant $C$ independent on $\eta$:
\begin{align*}
 &\Big| \E^{  { \mathbb{U}^{B}_H }}  \biggl( h(X^{\eta,  { \mathbb{U}^{B}_H }}_{1}) - \int_{0}^1  \tilde \lambda_*(X^{\eta,  { \mathbb{U}^{B}_H }}_{t},{\alpha}_{t}) dt  \biggr) -  \E^{  { \mathbb{U}^{B}_H }}  \biggl( h({X}^{0,  { \mathbb{U}^{B}_H }}_{1}) - \int_{0}^1  \tilde \lambda_*({X}^{0,  { \mathbb{U}^{B}_H }}_{t},{\alpha}_{t}) d t  \biggr)  \Big| \\ &\leq 
 C\,  \E^{  { \mathbb{U}^{B}_H }} \sup_{t \in [0,1]} | X^{ \eta,  { \mathbb{U}^{B}_H }} _t - {X}^{0,  { \mathbb{U}^{B}_H }} _t | .
\end{align*}

We have :
$$\begin{cases}
  d(X^{\eta,  { \mathbb{U}^{B}_H }} _t- {X}^{0,  { \mathbb{U}^{B}_H }}_t)= A(X^{\eta,  { \mathbb{U}^{B}_H }} _t- {X}^{0,  { \mathbb{U}^{B}_H }}_t) dt  -  [ R (X^{\eta,  { \mathbb{U}^{B}_H }} _t)- R(
  {X}^{0,\mathbb{U}^B}_t)] dW^1_t + \eta d B_t, & t\in [0,1] \\ \quad  X^{\eta,  { \mathbb{U}^{B}_H }} _0- {X}^{0,  { \mathbb{U}^{B}_H }}_0=0 &
\end{cases}$$
By standard estimates  based on the factorization method (see the proof of Theorem \ref{limite epsilon fisso}), we end up with:
\begin{equation*}
  \E^{  { \mathbb{U}^{B}_H }} \sup_{t \in [0,1]} | X^{ \eta,  { \mathbb{U}^{B}_H }} _t - {X}^{0,  { \mathbb{U}^{B}_H }} _t |^p \leq C_p \eta^p.
\end{equation*}
for all $p\geq 1 $ and a suitable constant $C_p$ independent of $\alpha$.
Thus  we can conclude that:
\begin{align*}
&
\Big|\inf_{   { \mathbb{U}^{B}_H }\in\,\mathbb{S}^{1,B}}\E^{\mathbb{U}^{B}_H }  \biggl( h(X^{\eta,  { \mathbb{U}^{B}_H }}_{1}) - \int_{0}^1  \tilde \lambda_*(X^{\eta,  { \mathbb{U}^{B}_H }}_{t},\alpha_{t}) dt  \biggr) -  \inf_{   { \mathbb{U}^{B}_H }\in\,\mathbb{S}^{1,B}} \E^{\mathbb{U}^{B}_H }  \biggl( h({X}^{0,\mathbb{U} ^{B,\mathfrak{ m }}}_{1} - \int_{0}^1  \tilde \lambda_*({X}^{0, { \mathbb{U}^{B}_H }}_{t},{\alpha}_{t}) dt  \biggr)\Big|   \\& \leq 
\sup_{   { \mathbb{U}^{B}_H }\in\,\mathbb{S}^{1,B}}  \Big| \E^{{ \mathbb{U}^{B}_H }}  \biggl( h({X}^{\eta,   { \mathbb{U}^{B}_H }}) - \int_{0}^1  \tilde \lambda_*({X}^{0,{ \mathbb{U}^{B}_H }}_{t},{\alpha}_{t}) dt  \biggr) -  \E^{\mathbb{U}^{B}_H }  \biggl( h({X}^{0,\mathbb{U}^{B}_H }_{1}) - \int_{0}^1  \tilde \lambda_*({X}^{0,{ \mathbb{U}^{B}_H }}_{t},{\alpha}_{t}) dt  \biggr)\Big|
\\&\leq  C \eta,
\end{align*} 
and the claim follows.  \finedim

\bigskip

 \begin{remark}
  In the special case in which the slow evolution is not perturbed  by the noise $(W^1)$ (equivalently $R\equiv 0$ in \eqref{lentovelocecontrollata}),  in Theorem \ref{main-caratt} the following characterization holds:
 $$ V(x_0) =  \inf_{ \mathbb{U}_H \in \, \mathbb{S}^0_H} \mathbb{E}^{\mathbb{U}^0_H}   \left( h({X}^{\mathbb{U}^0_H}_{1}) - \int_{0}^1  \tilde \lambda_*({X}^{\mathbb{U}^0_H}_s,{\alpha}_{s}^{\mathbb{U}^0_H}) ds \right) $$
 where $\mathbb{S}_H$ is the set of all $\mathbb{U}^0_H = (\Omega, \mathcal{F}, (\F_t), \mathbb{P}, (\alpha_t))$ where $\alpha$ is any $ (\F_t) $ progressively measurable process with values into the ball $B(0, M+1)\subset H$ of center 0 and radius $M+1$
and $ X^{\mathbb{U}^0_H}$ solves:
\begin{equation} \label{statolimite}d{X}_s= A{X}_s ds   + \alpha_s \, ds, \quad {X}_0=x_0.
\end{equation}
It is natural to think that the stochastic framework is here pleonastic and that the infimum above can be restricted to deterministic controls and trivial settings, namely:
\begin{equation}
V(x_0) = \inf_{ a}   \left( h({X}^{a}_{1}) - \int_{0}^1  \tilde \lambda_*({X}^{a}_s,a_s) ds \right)  
\end{equation}
where the above infimum is computed over all functions $a: [0,1] \rightarrow B(0,M+1)$ ($ B(0,M+1)$ being  the ball of $H$ centered in the origin of radius $M+1$) and $X^a$ is the mild solution of the deterministic evolution equation:
$$\frac{d}{dt}X_t= AX_t+a_t,\qquad X_0=x_0 $$
It is straight forward to see that 
$$\inf_{ \mathbb{U}^0_H \in \, \mathbb{S}_H} \mathbb{E}^{\mathbb{U}_H}   \left( h({X}^{\mathbb{U}^0_H}_{1}) - \int_{0}^1  \tilde \lambda_*({X}^{\mathbb{U}^0_H}_s,{\alpha}_{s}^{\mathbb{U}^0_H}) ds \right)\leq  \inf_{ a}   \left( h({X}^{a}_{1}) - \int_{0}^1  \tilde \lambda_*({X}^{a}_s,a_s) ds \right)  $$
To prove the converse let,  for any $\epsilon >0$, $ \mathbb{U}_H \in \mathbb{S}$ such that:
\begin{equation}
\mathbb{E}^{\mathbb{U}^{0}_H }   \left( h({X}^{\mathbb{U}^{0}_H }_{1}) - \int_{0}^1  \tilde \lambda_*({X}^{\mathbb{U}^{0}_H }_s,{\alpha}_{s}^{\mathbb{U}^{0}_H }) ds \right) \leq    \inf_{ \mathbb{U}^{0}_H  \in \, \mathbb{S}_H} \mathbb{E}^{\mathbb{U}^{0}_H }   \left( h({X}^{\mathbb{U}^{0}_H }_{1}) - \int_{0}^1  \tilde \lambda_*({X}^{\mathbb{U}^{0}_H }_s,{\alpha}_{s}^{\mathbb{U}^{0}_H }) ds \right)  + \epsilon.
\end{equation}
Then we have, recalling that equation \eqref{statolimite} can be solved pathwise:
\[  
\mathbb{P}^{\mathbb{U}^{0}_H } \left( \omega \in \Omega ^{ \mathbb{U}^{0}_H }: 
h({X}^{\mathbb{U}^{0}_H }_{1} (\omega)) - \int_{0}^1  \tilde \lambda_*({X}^{\mathbb{U}^{0}_H }_s (\omega),{\alpha}_{s}^{\mathbb{U}^{0}_H }(\omega)) ds
 \leq V(x_0)  + \epsilon\right)  >0
\]
To conclude is now enough to select $ \bar{\omega}$ in the above set and choose $ \bar{a}= \alpha ^{\mathbb{U}^{0}_H }(\bar{\omega})$. 
 \end{remark}

\begin{remark}\label{constrainedBSDEs}

Taking advantage of our main result  (see Theorem \ref{main-caratt}) we can further represent the singular limit $V(x_0)$ as the value at time $0$ of the minimal solution of a  BSDE with constrints on the martingale term (see \cite{KharroubiPham} for the definition for  and \cite{CossoGuatteriTessitore} for the infinite dimensional case). The bridge is given by the results in \cite{CossoGuatteriTessitore} allowing to represent the value function of a control problem by such a constrained   BSDE  without using viscosity solutions of the related HJB equation.  
 
\medskip 
 
First of all, to fit the framework in \cite{CossoGuatteriTessitore}  we notice that the control problem introduced in Theorem \ref{main-caratt} can be rewritten considering an unbounded set of controls. This is readily (and rather obviously) done by introduciong the class $\widehat{\mathbb{S}} ^{1}_H$  of $6$-uples  $\mathbb{U}_H=(\Omega,\F, (\F_t), \mathbb{P}, (W^1_t), (u_t))$ identical to the ones in $\mathbb{S} ^{1}_H$ with the only difference that here $u$ are not required to be bounded by $M+1$  and noticing that:
\begin{multline}
    {V}(x_0)= \inf_{\mathbb{U}_H \in \mathbb{S} ^ {1} _H} \mathbb{E}^{\mathbb{U}_H}  \left( h({X}^{\mathbb{U}_H}_{1}) - \int_{0}^1  \tilde \lambda_*({X}^{\mathbb{U}_H}_s,{\alpha}_{s}) ds \right) = \inf_{\widehat{\mathbb{U}}_H \in \widehat{\mathbb{S}}^ {1}_H} \mathbb{E}^{\mathbb{U}_H}      \left( h({X}^{\widehat{\mathbb{U}}_H}_1) - \int_{0}^1  \tilde \lambda_*({X}^{\widehat{\mathbb{U}}_H}_s,\Gamma(u_{s})) ds \right)
\end{multline}
where  $\Gamma(h):= \frac{h}{||h||} \min(||h||, M+1)$ for $h\in H$ and $X^{\widehat{\mathbb{U}}_H} $ solves:
$$d{X}_s= A{X}_s ds  - \, \Gamma({u}_s) ds  +  R ({X}_s) dW^1_s, \quad {X}_0=x_0.$$

Then we just have to remind the construction and the results in  \cite{CossoGuatteriTessitore}.

Given $x_0 \in H $ we consider the following system of forward-backward stochastic differential equations:
\begin{equation}\label{RandomizedFBSDE}
\begin{cases}
\mathcal{X} ^{x_0}_t\ =\displaystyle \ e^{tA} x_0+ \int_0^t   e^{(t-s)A} \Gamma(  S\widehat{\mathcal{W}}_s) ds +\, \int_0^t e^{(t-s)A}  R (\mathcal{X}^{x_0}_s) d\mathcal{W}_s, 
\\
  \mathcal{Y}^{x_0}_t=  h(\mathcal{X}^{x_0}_{1})  - \displaystyle\int_t^1 \tilde{\lambda}_* (\mathcal{X}^{x_0}_{s}, \Gamma( S\widehat{\mathcal{W}}_s))\, ds - K^{x_0}_1 + K^{x_0}_t- \int_t^1 \mathcal{Z}^{x_0}_s  \, d \mathcal{W}_s \ ,
    \end{cases}
\end{equation}
where  $ S : H \to H$ is an arbitrary trace class and injective linear operator with dense image and $ \mathcal{W}$ and $\widehat{\mathcal{W}}$ are two independent cylindrical Wiener processes with values in $H$ defined on a probability space satisfying the usual conditions. We denote by $(\F^0_t)$ the natural filtration of $(\mathcal{W}_t, \widehat{\mathcal{W}}_t)$ augmented.
Notice that, besides the two typical terms in the backward component,  the unknown  $K$ appears. Such process belongs to the set of real-valued $(\mathcal{F}^0_t)$- adapted nondecreasing continuous processes $K$ on $[0,T]$ such that $\E|K_T|^2< \infty$ and $K_0=0$.
\medskip

Then, see  \cite{CossoGuatteriTessitore} \S 4.2  the following holds:
\begin{itemize}
    \item  the forward equation in \eqref{RandomizedFBSDE} has a unique solution $(\mathcal{X}^{x_0})$ in  $L^{2}_{\F^{0}}(\Omega;{C} ([0,T];{H}))  $.
\item the backward equation in 
system  \eqref{RandomizedFBSDE} has a maximal solution $(\mathcal{Y}^{x_0}, \mathcal{Z}^{x_0}, K^{x_0})$ belonging to  the space  $L^{2}_{\F^{0}}(\Omega;{C} ([0,T];\mathbb{R}))\times L^{2}_{\F^{0}}(\Omega\times [0,T];\Xi^*)\times \mathcal{K}^2(0,T)$, maximal  in the sense that if there exists another solution $(\mathcal{Y}',\mathcal{Z}',{K}')$ belonging to the same functional spaces then
$\mathcal{Y}_t^{x_0} \geq \mathcal{Y}'_t$ for all $t\in [0,1]$, $\mathbb{P}$-a.s.
\item the following characterization of the singular limit $V(x_0)$  holds
\begin{equation}
   \mathcal{Y}^{x}_0  =  {V} (x_0).
\end{equation}
\end{itemize}
\end{remark}


\bibliographystyle{plain}
\bibliography{biblionew}

\end{document}